\documentclass[journal]{IEEEtran}

\usepackage{cite}
\usepackage{color}

\ifCLASSINFOpdf
  \usepackage[pdftex]{graphicx}
  \graphicspath{{../pdf/}{../jpeg/}}
\else
\fi

\usepackage{amsmath}
\usepackage{amsfonts}
\usepackage{amssymb}
\usepackage{cases}
\newtheorem{lemma}{Lemma}
\newtheorem{theorem}{Theorem}

\newtheorem{condition}{Condition}
\newtheorem{assumption}{Assumption}
\newtheorem{remark}{\textbf{Remark}}

\usepackage{multirow}
\usepackage{threeparttable}
\usepackage{float}

\usepackage{url}

\allowdisplaybreaks[4]

\usepackage{enumerate}
\usepackage{color}

\usepackage{comment}
\usepackage{ifthen}
\newboolean{showcomments}
\setboolean{showcomments}{true}
\newcommand{\re}[1]{\ifthenelse{\boolean{showcomments}}
	{ \textcolor{cyan}{(Re:  #1)}}{}}
 \newcommand{\fliu}[1]{\ifthenelse{\boolean{showcomments}}
 	{ \textcolor{green}{(FL:  #1)}}{}}
\newcommand{\slow}[1]{\ifthenelse{\boolean{showcomments}}
  { \textcolor{blue}{(SL:  #1)}}{}}

\begin{document}

\title{Optimal Power Flow in Stand-alone DC Microgrids}

\author{
		Jia~Li,
        Feng~Liu,~\IEEEmembership{Member,~IEEE,}
        Zhaojian~Wang,
        Steven~H~Low,~\IEEEmembership{Fellow,~IEEE,}
        Shengwei~Mei,~\IEEEmembership{Fellow,~IEEE,}
\thanks{
	
J. Li, F. Liu, Z. Wang and S. Mei are with the State Key Laboratory of Power System, the Department of Electrical Engineering, Tsinghua University, Beijing 100084, China.\{lfeng@tsinghua.edu.cn\}

S. H. Low is with the Department of Electrical Engineering, California
Institute of Technology, Pasadena, CA, USA, 91105. (email:slow@caltech.edu)}
}

\maketitle

\begin{abstract}
Direct-current microgrids (DC-MGs) can operate in either grid-connected or stand-alone mode. In particular, stand-alone DC-MG has many distinct applications. However, the optimal power flow problem of a stand-alone DC-MG is inherently non-convex. In this paper, the optimal power flow (OPF) problem of DC-MG is investigated considering convex relaxation based on second-order cone programming (SOCP). Mild assumptions are proposed to guarantee the exactness of relaxation, which only require uniform nodal voltage upper bounds and positive network loss. Furthermore, it is revealed that the exactness of SOCP relaxation of DC-MGs does not rely on either topology or operating mode of DC-MGs, and an optimal solution must be unique if it exists. If line constraints are considered, the exactness of SOCP relaxation may not hold. In this regard, two heuristic methods are proposed to give approximate solutions.
Simulations are conducted to confirm the theoretic results.
\end{abstract}


\IEEEpeerreviewmaketitle

\section{Introduction}
\IEEEPARstart{C}{ompared} with AC microgrids, Direct-Current microgrids (DC-MGs) have been recognized as an attractive alternative for numerous applications due to higher efficiency, more natural interface to many types of renewable energy resources and energy storage systems, better compliance with consumer electronics, etc. \cite{Dragicevic2014}. Additionally, when components are coupled around a DC bus, there are no issues with reactive power flow  and frequency stability, resulting in a notably less complex power system \cite{Dragicevic2016}. DC-MGs may operate in either grid-connected or stand-alone mode. 
The latter has many distinct applications in shipboard \cite{Reed2012,Jin2016a}, aircraft \cite{Magne2012,Magne2013}, automotive \cite{Emadi2006}, as well as the electricity supply of remote rural areas. In this regard, it is of great importance to investigate the optimal power flow of DC-MGs. 

Second-order cone programming (SOCP) has been extensively used in AC networks for solving the optimal power flow (OPF) problem \cite{Jabr2006,Sojoudi2012,Farivar2013a,Low2014,Low2013,Huang2016}. In \cite{Jabr2006}, radial distribution load flow is formulated as a conic quadratic optimization problem, and solved efficiently using interior-point methods. In \cite{Farivar2013a}, a two-step SOCP relaxation approach is proposed, which consist of angle relaxation and conic relaxation. The exactness of angle relaxation requires the so-called cyclic condition that the sum of angle differences on each cycle must be zero. With the cyclic condition satisfied,  the conic relaxation is exact, provided there are no upper bounds on loads. The formulations of OPF problem and their relaxations are summarized in \cite{Low2014}, and the sufficient conditions under which the convex relaxations are exact are presented in \cite{Low2013}. Another sufficient condition for the exact SOCP relaxation of AC OPF in radial distribution networks is proposed in \cite{Huang2016}, which requires that the allowed reverse power flow is only reactive or active, or none.

As the OPF problem of DC power network is inherently nonlinear and non-convex, SOCP is extended to realize the convexification of OPF problem in DC distributed networks \cite{Gan2014}. Such a result, however, relies on the assumption that there exists a large substation with unlimited capability to provide power injection and keep its nodal voltage constant. Unfortunately, a \emph{stand-alone} DC-MG is not the case because: 
	1) there is no substation with an unconstrained power injection, and every distributed generation therein has a limited capacity;
	2) there is no substation with a fixed nodal voltage, and all nodal voltages are allowed to vary within a certain range.  	
Such differences motivate us to extend the previous work \cite{Gan2014} to give rise to a better understanding of DC-MG in different operating modes.
To this end, the following steps are taken.

\begin{itemize}
	\item \textbf{Step 1}: \emph{ Equivalent Transformation}. By introducing slack variables, the original problem OPF1 is transformed equivalently into OPF2 with non-convex rank constraints. 
	\item \textbf{Step 2}: \emph{Seocnd-Order Conic Relaxation.} By removing the rank constraints, the non-convex problem OPF2 is relaxed into a convex SOCP problem, i.e., RL1. Moreover, extending the results of \cite{Gan2014}, this paper proves the exactness of SOCP relaxation under mild assumptions, which only require uniform nodal voltage upper bounds and positive network loss.
	\item \textbf{Step 3:} \emph{Equivalent Conversion.} By introducing alternative variables, RL1 is transformed equivalently into a branch flow model, namely RLS1, for the sake of improving numerical stability.
\end{itemize}

\begin{figure}[htb]
\centering
\includegraphics[width=3.5in]{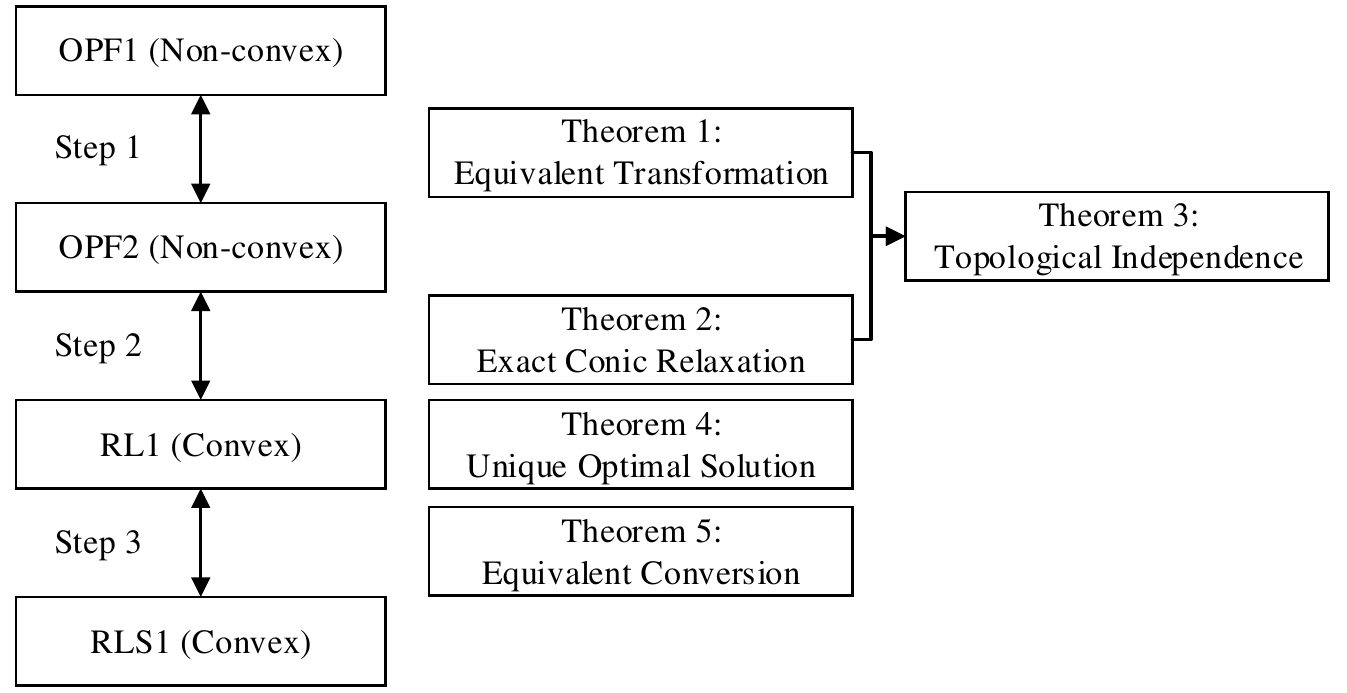}
\caption{Proposed solution approach.}
\label{fig:flow}
\end{figure}

The solution approach and the relation between the associated theoretic results are illustrated in Fig. \ref{fig:flow}.

Under the proposed mild assumptions, the exactness of SOCP relaxation and the uniqueness of optimal solution are revealed, which are independent of topologies and operating modes of DC-MGs. These properties of DC-MG are very helpful for optimal control and market design. When line constraints are considered, the situation turns to be much more complicated, since the rank constraint cannot be guaranteed. In this context, sufficient conditions are explored to justify the exactness of SOCP relaxation and two heuristic methods are suggested to construct feasible solutions when the rank constraint is not satisfied.

The rest of this paper is organized as follows. The OPF problem of a stand-alone DC-MG is formulated and equivalently transformed in Section \ref{sec:formulation}. The SOCP relaxation is given in Section \ref{sec:conic_relax}. Numerical studies are provided in Section \ref{sec:case}. The influence of line constraints are discussed in Section \ref{sec:discussion}. Section \ref{sec:conclusion} draws the conclusions.

\section{OPF Problem of DC-MGs}\label{sec:formulation}

\subsection{Basic Formulation}\label{sec:formulation1}
\begin{figure}[htb]
\centering
\includegraphics[width=2in]{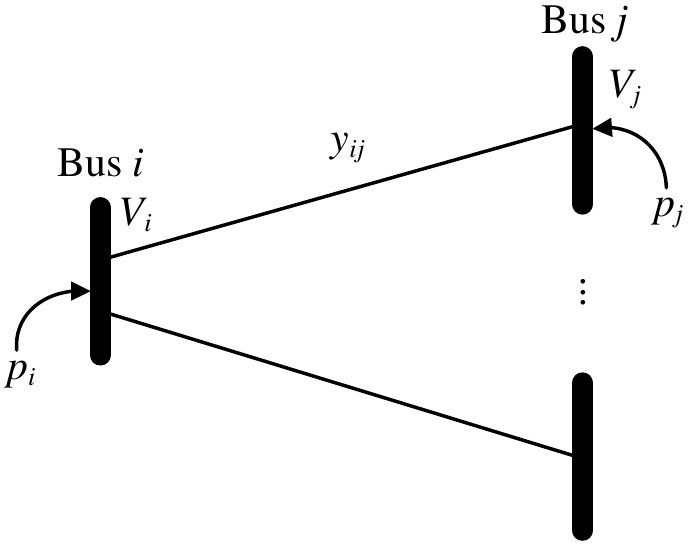}
\caption{Summary of notations.}
\label{fig:notation}
\end{figure}

Consider a graph $\mathcal{G}:=(\mathcal{N},\mathcal{E})$, where  $\mathcal{N}:=\{1,\cdots,n\}$ denotes the set of all buses and $\mathcal{E}$ denotes the set of all lines in the network of a DC-MG. $\mathcal{G}$ is assumed to be connected. Index the buses by $1,\cdots,n$ and abbreviate $\{i,j\}\in\mathcal{E}$ as $i\sim j$. Denote ($i\sim j$ \& $i<j$) by  $i\to j$. For each bus  $i\in \mathcal{N}$, denote $V_i$ as its voltage, and $p_i$ as its power injection. For each line $i\sim j $, $y_{ij}$ denotes its conductance. A letter without subscripts denotes a vector of the corresponding quantities, e.g., $V=[V_i]_{i\in \mathcal{N}}$. The notations are summarized in Fig.\ref{fig:notation}.
In a DC-MG, $V_i$, $p_i$, $y_{ij}$ are real numbers, and $y_{ij}>0$. Then the OPF problem of a stand-alone DC-MG reads
\begin{subequations}\label{eq:opf1}
\begin{align}
  \text{OPF1:} ~ \min ~ & h(p) = \sum_{i\in\mathcal{N}}f_i(p_i) \nonumber \\
  \text{over:} ~ & p,V; \nonumber \\
  \text{s.t.} ~
  & p_i = \sum_{j:j\sim i}V_i\left(V_i - V_j\right)y_{ij}, ~ i\in \mathcal{N}; \label{eq:opf1_p_i} \\
  & \underline{p}_i \le p_i \le \overline{p}_i, ~ i\in \mathcal{N}; \label{eq:opf1_p_gen} \\
  & \underline{V}_i \le V_i \le \overline{V}_i, ~ i\in  \mathcal{N}. \label{eq:opf1_v}
\end{align}
\end{subequations}
Here, $f_i(p_i)$ is strictly increasing in $p_i$. Eq. \eqref{eq:opf1_p_i} is the power injection equation for bus $i$. The nodal voltages are constrained by \eqref{eq:opf1_v} with the lower bound  $\underline{V}_i>0$ and upper bound $\overline{V}_i$. The power injections are constrained by \eqref{eq:opf1_p_gen} with the lower bound $\underline{p}_i$ and upper bound $\overline{p}_i$. It is assumed that $\underline{p}_i\le 0$, which indicates the following two cases:

\begin{enumerate}

\item  Bus $i$ is a pure generation bus without load. In this case, $\underline{p}_i= 0$ as the generators in a DC-MG  can be turned off; 

\item  Bus $i$ is a pure load bus without generation or a mixture power injection of both load and generation. In this case, there is $\underline{p}_i< 0$.

\end{enumerate}

In fact, the model of a grid-connected DC-MG \cite{Gan2014} can be viewed as a special case of OPF1 including a substation bus with an unconstrained power injection ($\underline{p}_0=-\infty,\overline{p}_0=\infty$) and a fixed nodal voltage ($V_0=V_0^{\text{ref}}$).

Note that, in a DC-MG, line constraints usually do not bind in the normal operation condition, since the network is overprovisioned. Therefore, in this paper, we first consider the OPF problem of DC-MGs without line constraints. Then, we discuss the influence of line constraints on the exactness of SOCP relaxation. Throughout the paper, we do not assume any specific topology of the power network.

\subsection{Equivalent Transformation}\label{sec:eq_trans}
The proposed OPF problem \eqref{eq:opf1} is a non-linear non-convex problem. By introducing slack variables, it can be transformed into an equivalent counterpart, where the non-convex power injection equation \eqref{eq:opf1_p_i} is converted into a rank constraint.

Introduce slack variables to formulate a map $f$ such that
\begin{subequations}\label{eq:slack_var}
\begin{numcases}
  {f:=}
  v_i = V_i^2, & $i\in \mathcal{N}$; \label{eq:slack_v} \\
  W_{ij} = V_iV_j, & $i\sim j$; \label{eq:slack_w}
\end{numcases}
\end{subequations}
and define a matrix
\begin{equation}
R_{ij}:=\begin{bmatrix}
v_i & W_{ij} \\
W_{ji} & v_j
\end{bmatrix}
\end{equation}\label{eq:trans}
for every $i\to j$. 
Then OPF1 \eqref{eq:opf1} is transformed into 
\begin{subequations}\label{eq:opf2}
\begin{align}
  \text{OPF2:} ~ \min ~ & h(p) \nonumber \\
  \text{over:} ~ & p,v,W \nonumber \\
  \text{s.t.} ~ & p_i = \sum_{j:j\sim i}\left(v_i - W_{ij}\right)y_{ij}, ~ i\in \mathcal{N} \label{eq:opf2_p_i} \\
   & \underline{p}_i \le p_i \le \overline{p}_i, ~ i\in \mathcal{N} \label{eq:opf2_p_gen} \\
  & \underline{V}_i^2 \le v_i \le \overline{V}_i^2, ~ i\in  \mathcal{N} \label{eq:opf2_v} \\
  & W_{ij} \ge 0, ~ i\to j \label{eq:opf2_w_0} \\
  & W_{ij} = W_{ji}, ~ i\to j \label{eq:opf2_w} \\
  & R_{ij} \succeq 0, ~ i\to j \label{eq:opf2_psd} \\
  & \text{rank}(R_{ij})=1, ~ i\to j  \label{eq:opf2_rank}
\end{align}
\end{subequations}
where the non-convexity in \eqref{eq:opf1_p_i} (in OPF1) is converted into the non-convexity in the rank constraint \eqref{eq:opf2_rank} (in OPF2). $R_{ij}$ is positive semidefinite as shown in \eqref{eq:opf2_psd}. 

\begin{theorem}\label{theorem:eq_trans}
OPF1 and OPF2 are equivalent.
\end{theorem}

To prove Theorem \ref{theorem:eq_trans}, we first give the following lemma.
\begin{lemma}\label{lemma:trans}
Given $v_i>0$ for $i\in \mathcal{N}$ and $W_{ij}\ge 0$ for $i\to j$, let  $W_{ij}=W_{ji}$ for $i\to j$. If $\text{rank}(R_{ij})=1$ for $i\to j$, then there exists a unique $V$ satisfying $V_i>0$ for $i\in \mathcal{N}$ and \eqref{eq:slack_var}. Moreover, $V$ is  determined by $V_i = \sqrt{v_i}$ for $i \in \mathcal{N}$.
\end{lemma}

The proof of Lemma \ref{lemma:trans} can be found in Appendix \ref{sec:proof_lemma_trans}. 

Lemma \ref{lemma:trans} implies that {for each $(v,W)$, there exists a unique $V$ that satisfies the map given by \eqref{eq:slack_var}}.
Differing from \cite{Gan2014}, we do not have a substation bus with a fixed voltage. Instead, we assume $V_i>0$ for all $i\in\mathcal{N}$.  
With Lemma \ref{lemma:trans}, Theorem \ref{theorem:eq_trans} can be proven. The proof can be found in  Appendix \ref{sec:proof_thm_trans}. 

It is worth noting that Theorem \ref{theorem:eq_trans} does not rely on the topology of network.

\section{Exactness of Conic Relaxation}\label{sec:conic_relax}
\subsection{SOCP Relaxation of OPF in DC-MGs}\label{sec:RL1}
By removing \eqref{eq:opf2_rank}, the original non-convex problem OPF2 is transformed to an SOCP problem (named as RL1) as below.
\begin{align}
  \text{RL1:} ~ \min ~ & h(p) \nonumber \\
  \text{over:} ~ & p,v,W \nonumber \\
  \text{s.t.} ~ & \eqref{eq:opf2_p_i}\text{ -- }\eqref{eq:opf2_psd}\nonumber
\end{align}

The only difference between RL1 and OPF2 is that RL1 has no constraint \eqref{eq:opf2_rank}. Therefore, RL1 is \emph{exact}, provided that its every optimal solution satisfies \eqref{eq:opf2_rank}. To ensure the exactness of  conic relaxation, additional assumptions are required.

\subsection{Assumptions}\label{sec:condition}

Throughout the paper, we make the following assumptions.
\begin{assumption}\label{assumption1}
  $\overline{V}_1=\overline{V}_2=\cdots\overline{V}_n > 0$.
\end{assumption}

\begin{assumption}\label{assumption2}
 $\sum_{i\in\mathcal{N}}p_i>0$.
\end{assumption}

Assumption \ref{assumption1} requires all the nodal voltages have the same upper bounds, which is reasonable in DC-MGs, since the scale of system is usually small. Assumption \ref{assumption2} is trivial as it means the total network loss is positive.
Such assumptions relax those in \cite{Gan2014}, which require negative power injection lower bounds and an unconstrained power injection, to admit the features of stand-alone DC-MGs.

\subsection{Exactness of Conic Relaxtion}\label{sec:theory}
With Assumptions \ref{assumption1} and \ref{assumption2} mentioned above, we have the following main theorem:

\begin{theorem}\label{theorem:exact_relax}
   RL1 is  exact if Assumptions \ref{assumption1} and \ref{assumption2} hold for OPF1 (equivalently OPF2). 
\end{theorem}

Theorem \ref{theorem:exact_relax} claims that RL1 is an exact SOCP relaxation of OPF2 (equivalently OPF1) under Assumptions \ref{assumption1} and \ref{assumption2}. To prove this theorem, we introduce the following lemmas.

\begin{lemma}\label{lemma:v}
	Assume Assumption \ref{assumption1} holds. Let $(p,v,W)$ be feasible for RL1 but violate the rank constraint \eqref{eq:opf2_rank} on a certain line $(s\to t)\in\mathcal{E}$. If $p_s=\underline{p}_s$, then $v_s<\overline{V}_s^2$. Meanwhile, if $p_t=\underline{p}_t$, then $v_t<\overline{V}_t^2$.
\end{lemma}

The proof of Lemma \ref{lemma:v} can be found in Appendix \ref{sec:proof_lemma_v}. 

Lemma \ref{lemma:v} implies that, for each bus, the power injection's lower bound and the nodal voltage's upper bound cannot bind at the same time. Therefore, for any feasible solution $(p,v,W)$ of RL1, if it violates the rank constraint \eqref{eq:opf2_rank} for a certain line $(s\to t)\in\mathcal{E}$ and the constraint $p_s\ge\underline{p}_s$ (or $p_t\ge\underline{p}_t$) is binding, then $v_s\le \overline{V}_s^2$ (or $v_t\le \overline{V}_t^2$) cannot bind. 

\begin{lemma}\label{lemma:le}
Assume Assumptions \ref{assumption1} and \ref{assumption2} hold for OPF1  and let $(p,v,W)$ be a feasible solution to RL1. If
\begin{enumerate}
  \item $(p,v,W)$ violates the rank constraint \eqref{eq:opf2_rank} on a certain line $(s\to t)\in\mathcal{E}$;
  \item $p_s>\underline{p}_s$, $p_t>\underline{p}_t$;
\end{enumerate}
then there  exists another feasible solution  $(p',v,W')$ that
\begin{enumerate}
  \item satisfies \eqref{eq:opf2_p_i} -- \eqref{eq:opf2_psd};
  \item satisfies $h(p')<h(p)$.
\end{enumerate}
\end{lemma}

The proof of Lemma \ref{lemma:le} can be found in Appendix \ref {sec:proof_lemma_le}. 

Lemma \ref{lemma:le} says that if a feasible point $(p,v,W)$ violates the rank constraint \eqref{eq:opf2_rank} for a certain line $s\to t$, while both $p_s$ and $p_t$ are not binding, then we can always find another feasible point $(p',v,W')$ with a better objective value. It implies that if the optimal solution $(p^*,v^*,W^*)$ to RL1 violates the rank constraint \eqref{eq:opf2_rank} for a certain line $s\to t$, then at least one of $p^*_s$ and $p^*_t$ must have reached its lower bound.

\begin{lemma} \label{lemma:le_eq}
Assume Assumptions \ref{assumption1} and \ref{assumption2} hold for OPF1 and let $(p,v,W)$ be a feasible solution to RL1. If
\begin{enumerate}
  \item $(p,v,W)$ violates the rank constraint  \eqref{eq:opf2_rank} on a certain line $(s\to t)\in\mathcal{E}$;
  \item either ($p_s=\underline{p}_s$ \& $p_t>\underline{p}_t$) or ($p_s>\underline{p}_s$ \& $p_t=\underline{p}_t$);
\end{enumerate}
then there exists another feasible solution $(p',v',W')$ that
\begin{enumerate}
  \item satisfies \eqref{eq:opf2_p_i} -- \eqref{eq:opf2_psd};
  \item satisfies $h(p')<h(p)$.
\end{enumerate}
\end{lemma}

The proof of Lemma \ref{lemma:le_eq} can be found in Appendix \ref{sec:proof_lemma_le_eq}. 

Lemma \ref{lemma:le_eq} means that if a feasible point $(p,v,W)$ violates the rank constraint \eqref{eq:opf2_rank} for a certain line $s\to t$, while either $p_s$ or $p_t$ is binding, then we can always find another feasible point $(p',v',W')$ with a better objective value. 
Lemma \ref{lemma:le} and \ref{lemma:le_eq} imply that if the optimal solution $(p^*,v^*,W^*)$ to RL1 violates the rank constraint \eqref{eq:opf2_rank} for a certain line $s\to t$, then it must satisfy $p^*_s=\underline{p}_s$ and $p^*_t=\underline{p}_t$. Otherwise, we can always find a better solution.

\begin{lemma}\label{lemma:eq}
Assume Assumptions \ref{assumption1} and \ref{assumption2} hold for OPF1 and let $(p,v,W)$ be a feasible solution to RL1. If
\begin{enumerate}
  \item $(p,v,W)$ violates the rank constraint \eqref{eq:opf2_rank} on a certain line $(s\to t)\in\mathcal{E}$;
  \item $p_s=\underline{p}_s$ and $p_t=\underline{p}_t$;
\end{enumerate}
then there always exists $(p,v',W')$ that
\begin{enumerate}
  \item satisfies \eqref{eq:opf2_p_i} -- \eqref{eq:opf2_psd};
  \item violates rank constraint \eqref{eq:opf2_rank} for all the neighbouring lines of $s\to t$, i.e. $i\to j$ with $\{\{i\},\{j\}\}\cap\{\{s\},\{t\}\}\ne \emptyset$.
\end{enumerate}
\end{lemma}

The proof of Lemma \ref{lemma:eq} can be found in Appendix \ref {sec:proof_lemma_eq}.

Lemma \ref{lemma:eq} says that if a feasible point $(p,v,W)$ violates the rank constraint \eqref{eq:opf2_rank} for a certain line $s\to t$, while the corresponding power injections' lower bounds are binding, we can always find another feasible point $(p,v',W')$ which  violates the rank constraint \eqref{eq:opf2_rank} for $s\to t$ and all its neighboring lines. It should be noted that in Lemma \ref{lemma:eq}, the construction of feasible point $(p,v',W')$ does not change $p$. Since the network $\mathcal{G}$ is connected, Lemma \ref{lemma:eq} implies that we can continue such propagation to obtain a feasible point that violates the rank constraint \eqref{eq:opf2_rank} for all the lines, without changing $p$.

As mentioned before, if the optimal solution $(p^*,v^*,W^*)$ to RL1 violates the rank constraint \eqref{eq:opf2_rank} for a certain line $s\to t$, Lemma \ref{lemma:le} and \ref{lemma:le_eq} ensure that it must satisfy $p^*_s=\underline{p}_s$ and $p^*_t=\underline{p}_t$. In this case, according to Lemma \ref{lemma:eq}, we can always find a feasible $(p',v',W')$ which satisfies $p'=p^*$ and violates the rank constraint \eqref{eq:opf2_rank} for all the lines. Following this idea, we can prove Theorem \ref{theorem:exact_relax} based on Lemma \ref{lemma:v} -- \ref{lemma:eq}. 

The proof of Theorem \ref{theorem:exact_relax} can be found in Appendix \ref{proof_thm_relax}.

\subsection{Topological Independence}
Recall the SOCP relaxation in AC networks \cite{Farivar2013a}. The approach consists of two relaxation steps: angle relaxation and conic relaxation. Similarly, our method also contains two steps: equivalent transformation and conic relaxation. Differing  from AC networks, DC networks do not involve voltage angles. The first step in our method (i.e., equivalent transformation) is exact, as Theorem \ref{theorem:eq_trans} states. However, in the second step, directly removing the rank constraints may result in inexactness. Thus, additional assumptions (Assumptions \ref{assumption1},  \ref{assumption2}) are made to ensure the exactness of  conic relaxation. By noting that none of the two steps depends on specific network topologies, we directly have the following theorem:

\begin{theorem}\label{theorem:topo_ind}
Assume Assumptions \ref{assumption1} and  \ref{assumption2} hold. Then the exactness of RL1 is independent of network topologies.
\end{theorem}

\begin{remark}
\emph{In terms of a grid-connected DC-MG, it has also been demonstrated that the SOCP relaxation is topology-independent \cite{Gan2014}. Actually, when a DC-MG works in grid-connected state, one can simply assign a substation bus in RL1 by letting $\underline{p}_0=-\infty,\overline{p}_0=+\infty$ and $v_0=[V_0^{\text{ref}}]^2$. Hence, combining the results in \cite{Gan2014} and Theorem \ref{theorem:topo_ind} immediately  concludes that the SOCP relaxation of the OPF problem of DC-MGs does not rely on the topology of network, whether the DC-MG works in grid-connected or stand-alone mode.}
\end{remark}

\subsection{Uniqueness of the Optimal Solution}
Next we show the  the optimal solution to RL1 is unique. 
\begin{theorem}
\label{theorem:RL1_unique}
Assume Assumptions \ref{assumption1} and  \ref{assumption2} hold for RL1. Then the optimal solution to RL1 is unique.
\end{theorem}

The proof of Theorem \ref{theorem:RL1_unique} can be found in Appendix \ref{proof_thm_unique}.

In a grid-connected DC-MG, it has also been proven that the SOCP problem has at most one optimal solution \cite{Gan2014}. Thus, the uniqueness of optimal solution to the SOCP problem does not depend on the operating mode of DC-MG.

\subsection{Branch Flow Model}
In an optimal solution to RL1, $v_i$ and $W_{ij}$ may be numerically close to each other, since the range of nodal voltage is small (usually 0.95$\sim$1.05 p.u.) and $\text{Rank}(R_{ij})=1$ is satisfied, which implies that $v_iv_j=W_{ij}W_{ji}$.
Thus, RL1 is ill-conditioned since equation \eqref{eq:opf2_p_i} requires the subtractions of $v_i$ and $W_{ij}$. However, such subtractions can be avoided by converting RL1 into a branch flow model, so that the numerical stability is improved. In light of \cite{Gan2014}, by defining $z_{ij}:=1/y_{ij}$ and adopting alternative variables $P$, $l$, RL1 can be converted into a branch flow model via the map $g:(v,W)\mapsto (v',P,l)$  defined as below.
\begin{subequations}\label{eq:convert}
\begin{numcases}
  {g:=}
  v'_i = v_i, & $i\in\mathcal{N}$; \\
  P_{ij} = \left(v_i - W_{ij}\right)y_{ij}, & $i\to j$; \\
  l_{ij} = y_{ij}^2\left(v_i - W_{ij} - W_{ji} + v_j\right), & $i\to j$.
\end{numcases}
\end{subequations}

With $g$, RL1 is converted into the following optimization problem with branch flow model (BFM):
\begin{subequations}\label{eq:rls1}
\begin{align}
  \text{RLS1:} ~ \min ~ & h(p) \nonumber \\
  \text{over:} ~ & p,P,v,l; \nonumber \\
  \text{s.t.} ~ & p_i = \sum_{j:j\sim i}P_{ij},~ i\in \mathcal{N}; \label{eq:rls1_p_i} \\
  & \underline{p}_i \le p_i \le \overline{p}_i, ~ i\in \mathcal{N}; \label{eq:rls1_p_gen} \\
  & \underline{V}_i^2 \le v_i \le \overline{V}_i^2, ~ i\in  \mathcal{N}; \label{eq:rls1_v} \\
  & P_{ij} + P_{ji} = z_{ij}l_{ij}, ~ i\to j; \\
  & v_i - v_j = z_{ij}\left(P_{ij} - P_{ji}\right), ~ i\to j; \label{eq:rls1_v_p} \\
  & l_{ij} \ge \frac{P_{ij}^2}{v_i}, ~ i\sim j \label{eq:rls1_rank}
\end{align}
\end{subequations}
where $P_{ij}$ denote the power flow through line $i\to j$, and $l_{ij}$ denote the magnitude square of the current through line $i\to j$.

\begin{theorem}\label{theorem:eq_conv}
  RL1 and RLS1 are equivalent.
\end{theorem}

The proof of Theorem \ref{theorem:eq_conv} can be found in Appendix \ref{proof_thm_eq_conv}.

Since RL1 is a convex relaxation of the non-convex OPF2 by removing the rank constraint \eqref{eq:opf2_rank}, Theorem \ref{theorem:eq_conv} implies that RLS1 is the convex relaxation of another non-convex problem, namely, OPF3, which is obtained by converting OPF2 using the one-to-one map $g$ \eqref{eq:convert}.
\begin{align}
  \text{OPF3:} ~ \min ~ & h(p) \nonumber \\
  \text{over:} ~ & p,P,v,l \nonumber \\
  \text{s.t.} ~ & \eqref{eq:rls1_p_i}\text{ -- }\eqref{eq:rls1_v_p}\nonumber \\
  & l_{ij} = \frac{P_{ij}^2}{v_i}, i\sim j \nonumber
\end{align}

\begin{figure}[htb]
\centering
\includegraphics[width=3.5in]{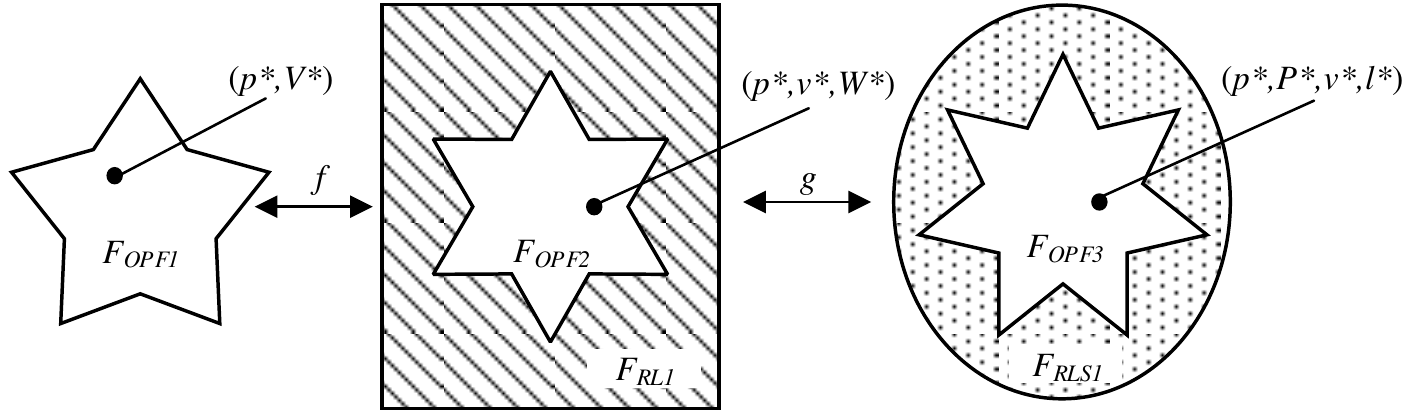}
\caption{Relationship between different feasible sets.}
\label{fig:set}
\end{figure}

Under Assumptions \ref{assumption1} and \ref{assumption2}, let $\mathcal{F}$ denote a feasible set of a certain optimization problem, which is indicated by the subscript.  
Then the relationship between different feasible sets is shown as below:
\begin{align*}
\mathcal{F}_{OPF2} &= f(\mathcal{F}_{OPF1}) \subset \mathcal{F}_{RL1} \\
\mathcal{F}_{RLS1} &= g(\mathcal{F}_{RL1})
\end{align*}
where $f$ is the one-to-one map given in Theorem \ref{theorem:eq_trans}, and $g$ is the one-to-one map given in Theorem \ref{theorem:eq_conv}. The relationship between different feasible sets is depicted in Fig. \ref{fig:set}. $\mathcal{F}_{OPF1}$ is a non-convex region and $(p^*,V^*)$ denotes the optimal solution to OPF1. $\mathcal{F}_{OPF1}$ is transformed equivalently into another non-convex region $\mathcal{F}_{OPF2}$ by the map $f$. $\mathcal{F}_{OPF2}$ is convexified by removing the rank constraint \eqref{eq:opf2_rank}, yielding $\mathcal{F}_{RL1}$, which is larger than $\mathcal{F}_{OPF2}$. The exactness of RL1 ensures that the optimal solution to RL1, i.e., $(p^*,v^*,W^*)$ is not in the shaded region of $\mathcal{F}_{RL1}$. Since the objective function of OPF2 and RL1 are the same, $(p^*,v^*,W^*)$ is also the optimal solution to OPF2. Applying the one-to-one map $f$, $(p^*,v^*,W^*)$ is transformed into $(p^*,V^*)$. The one-to-one map $g$ converts $\mathcal{F}_{RL1}$ into an equivalent convex region $\mathcal{F}_{RLS1}$, and converts $\mathcal{F}_{OPF2}$ into an equivalent non-convex region $\mathcal{F}_{OPF3}$.
$(p^*,P^*,v^*,l^*)$ is the optimal solution to RLS1, which is also the optimal solution to OPF3, since RLS1 is exact. Moreover, $(p^*,P^*,v^*,l^*)$ can be converted equivalently into $(p^*,v^*,W^*)$ using the map $g$.

Extending the work in \cite{Gan2014}, we have shown that under the proposed assumptions, the OPF problem in DC-MGs can be converted equivalently into a convex SOCP problem, regardless of topologies or operating modes. Thereby solving the relaxed convex SOCP can obtain the global optimum with theoretic guarantee.  Furthermore, we have shown that the SOCP problem has at most one optimal solution, whether the DC-MG is working in grid-connected or stand-alone mode.
\begin{remark}
\emph{	In OPF1$\sim$OPF3, line flow constraints are not considered. Whereas we have not found any violation of line constraints in all our tests,  we cannot theoretically exclude the possibility of such violation. In such a circumstance, approximate solutions can be heuristically constructed. We will discuss the situation this issue later on in Section \ref{sec:discussion}.  
}
\end{remark}

\section{Case Studies}\label{sec:case}

\subsection{16-bus System}
\begin{figure}[htb]
\centering
\includegraphics[width=3in]{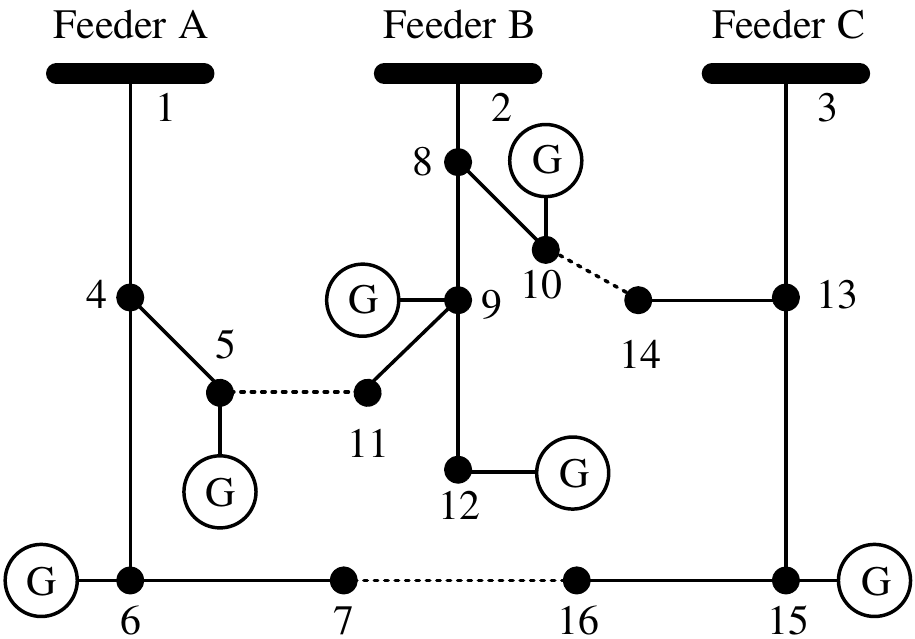}
\caption{16-bus system with DGs.}
\label{fig:case16}
\end{figure}

Tests are conducted on the modified 16-bus system \cite{Civanlar1988} to demonstrate the efficacy of the proposed method when the system works in different operating modes and network topologies. The three-feeder system is shown in Fig. \ref{fig:case16} where 6 distributed generators (DGs) are added with the same capacity of 5MW. The dashed lines represent tie lines with tie breakers. The system is able to work in two modes:
\begin{itemize}
\item \textbf{Grid-connected}: the system is connected to a power grid via all the three feeders (Feeders A, B and C).
\item \textbf{Stand-alone}: all the three feeders are disconnected from the grid.
\end{itemize} 
Additionally, the network can switch between tree and mesh topologies by opening or closing the tie breakers.
\begin{itemize}
\item \textbf{Tree topology}: all the tie breakers are switched off.
\item \textbf{Mesh topology}: all the tie breakers are switched on.
\end{itemize}

Thus, the following four cases are studied.

1) \textbf{GT}: grid-connected mode with a tree topology.

2) \textbf{GM}: grid-connected mode with a mesh topology.

3) \textbf{ST}: stand-alone mode with a tree topology.

4) \textbf{SM}: stand-alone mode with a mesh topology.

The voltage lower and upper bounds are set as 0.95 p.u. and 1.05 p.u., respectively. The objective is to minimize total network loss and the problems are solved using MOSEK. The results are listed in Table \ref{tab:case16}.

\begin{table}[htbp]
  \centering
  \caption{Results of 16-Bus System}
    \begin{tabular}{ccccc}
    \hline
          & GT    & GM    & ST    & SM \\
    \hline
    Exactness & 5.73E-9 & 1.46E-9 & 5.59E-10 & 1.72E-10 \\
    Objective Value (p.u.)   & 0.012 & 0.009 & 0.017 & 0.013 \\
    Computation Time (s) & 0.13 & 0.20 & 0.12 & 0.12 \\
    \hline
    \end{tabular}%
  \label{tab:case16}%
\end{table}%

At a numerical RLS1 solution $(p,v,W)$, $R_{ij}$ can be obtained for each $i\to j$. If RLS1 is exact, then for all $i\to j$, we should have $\text{rank}(R_{ij})=1$, i.e., the difference $D_{ij}:=v_iv_j - W_{ij}W_{ji}=0$. Hence the smaller $D_{ij}$ indicates, the closer $R_{ij}$ is to rank one. The row ``Exactness'' lists the maximum value of $D_{ij}$ for all $i\to j$, indicating RLS1 is exact for all the test networks. Moreover, it indicates that the proposed SOCP relaxation does not depend on the operating mode or network topology. The objective values indicate that in this system, mesh topology reduces network loss in both grid-connected and stand-alone modes, since mesh topology is able to support higher nodal voltages. Additionally, working in the same network topology, grid-connected mode experiences less network loss than stand-alone mode, since more power is injected from the feeders to support higher nodal voltages. The results of power injections and nodal voltages are listed in Table \ref{tab:case16_p_v}.

\begin{table}[htbp]
  \centering
\scriptsize
  \caption{Results of Power Injections and Nodal Voltages}
    \begin{tabular}{crrrrrrrr}
    \hline
          & \multicolumn{4}{c}{Power Injection (p.u.)} & \multicolumn{4}{c}{Nodal Voltage (p.u.)} \\
    Bus & GT    & GM    & ST    & SM    & GT    & GM    & ST    & SM \\
    \hline
    1     & 0.072 & 0.072 & 0 & 0 & 1.050 & 1.050 & 1.042 & 1.042 \\
    2     & 0.144 & 0.131 & 0 & 0 & 1.050 & 1.050 & 1.026 & 1.029 \\
    3     & 0.085 & 0.041 & 0 & 0 & 1.050 & 1.050 & 1.034 & 1.044 \\
    4     & -0.200 & -0.200 & -0.200 & -0.200 & 1.045 & 1.045 & 1.042 & 1.042 \\
    5     & 0.068 & 0.132 & 0.107 & 0.192 & 1.050 & 1.050 & 1.050 & 1.050 \\
    6     & 0.211 & 0.226 & 0.246 & 0.260 & 1.050 & 1.050 & 1.050 & 1.050 \\
    7     & -0.150 & -0.150 & -0.150 & -0.150 & 1.044 & 1.044 & 1.044 & 1.044 \\
    8     & -0.400 & -0.400 & -0.400 & -0.400 & 1.035 & 1.036 & 1.026 & 1.029 \\
    9     & 0.100 & 0.100 & 0.100 & 0.100 & 1.044 & 1.047 & 1.039 & 1.044 \\
    10    & 0.144 & 0.235 & 0.231 & 0.324 & 1.050 & 1.050 & 1.050 & 1.050 \\
    11    & -0.060 & -0.060 & -0.060 & -0.060 & 1.038 & 1.048 & 1.033 & 1.047 \\
    12    & 0.079 & 0.038 & 0.139 & 0.078 & 1.050 & 1.050 & 1.050 & 1.050 \\
    13    & -0.100 & -0.100 & -0.100 & -0.100 & 1.041 & 1.046 & 1.034 & 1.044 \\
    14    & -0.100 & -0.100 & -0.100 & -0.100 & 1.032 & 1.046 & 1.026 & 1.045 \\
    15    & 0.329 & 0.253 & 0.415 & 0.279 & 1.050 & 1.050 & 1.050 & 1.050 \\
    16    & -0.210 & -0.210 & -0.210 & -0.210 & 1.042 & 1.042 & 1.042 & 1.042 \\
    \hline
    \end{tabular}%
  \label{tab:case16_p_v}%
\end{table}%

\subsection{Exactness and Comparison}
More tests are conducted on both mesh and tree networks to check the exactness of the SOCP relaxation. In addition, the objective value and computation time of the relaxed model (i.e., RLS1) and the original model (i.e., OPF1) are compared to show the efficacy of the relaxation. The objective of both models is to minimize total network loss. The mesh networks \cite{Gan2014} are modified from MATPOWER by ignoring line reactances. All the line resistances are reduced to 10\% of the original values to simulate the DC-MG condition. Particularly, the zero resistances are reset as $10^{-3}$ p.u.. Additionally, the IEEE 118-bus system is applied to show the scalability of the proposed method. Four radial distribution networks in the literature are also used to verify the topological independence of the relaxation. In these systems, all the line resistances are also reduced to 10\% of the original values to simulate the DC-MG condition. In all the test systems, the voltage lower and upper bounds are set as 0.95 p.u. and 1.05 p.u., respectively. RLS1 is solved using MOSEK, while OPF1 is solved using IPOPT. The results are listed in Table \ref{tab:compare}.

\begin{table}[htbp]  
\centering
\begin{threeparttable}
\scriptsize
  \caption{Exactness of RLS1 and Comparison of Two Models}
    \begin{tabular}{rcrrrrr}
    \hline
          &       &  Exactness     & \multicolumn{2}{c}{Objective Value (p.u.)} & \multicolumn{2}{c}{Time (s)} \\
    Topology & System & of RLS1 & RLS1  & OPF1  & RLS1  & OPF1 \\
    \hline
    \multicolumn{1}{c}{\multirow{5}[0]{*}{Mesh}} & case6ww & 1.24E-10 & 3.17E-03 & 3.17E-03 & 0.14  & 0.13 \\
    \multicolumn{1}{c}{} & case9 & 7.17E-12 & 5.72E-03 & 5.72E-03 & 0.11  & 0.14 \\
    \multicolumn{1}{c}{} & case\_ieee30 & 2.37E-11 & 1.52E-03 & 1.52E-03 & 0.14  & 0.19 \\
    \multicolumn{1}{c}{} & case39 & 3.64E-11 & 1.30E-01 & 1.30E-01 & 0.16  & 0.16 \\
    \multicolumn{1}{c}{} & case118 & 6.38E-11 & 7.98E-03 & 7.98E-03 & 0.21  & 0.23 \\
    \hline
    \multicolumn{1}{c}{\multirow{3}[0]{*}{Tree}} & 33-bus \cite{Baran1989} \tnote{a} & 1.28E-11 & 1.47E-01 & 1.47E-01 & 0.13  & 0.14 \\
    \multicolumn{1}{c}{} & 70-bus \cite{Das2006} \tnote{b} & 5.35E-12 & 1.10E-01 & 1.10E-01 & 0.23  & 0.21 \\
    \multicolumn{1}{c}{} & 94-bus \cite{Ching-TzongSu2003} \tnote{c} & 3.09E-11 & 3.63E-01 & 3.63E-01 & 0.22  & 0.36 \\
    \hline
    \end{tabular}%
  \label{tab:compare}%
      \begin{tablenotes}
        \item [a] 6 DGs are added at bus 5, 10, 15, 20, 25, 30, with the capacity of 50kW.
        \item [b] 13 DGs are added at bus 5, 10, 15, 20, 25, 30, 35, 40, 45, 50, 55, 60, 65, with the capacity of 50kW.
        \item [c] 13 DGs are added at bus 5, 10, 15, 20, 25, 30, 35, 40, 45, 50, 55, 60, 65, with the capacity of 50kW.
    \end{tablenotes}
\end{threeparttable} 
\end{table}%

It is observed that RLS1 is exact for both mesh and tree topologies. The scale of a DC-MG is small, however, more and more DC-MGs may be built in the future to integrate distributed renewable energy generation and electric vehicles. The results of computation time implies that the proposed method may be applied in large scale systems, for example, the cluster of DC-MGs.

We also tested the model with line constraints (i.e., RLS2), on the above systems. However, we have not found any case where the optimal solution violates the rank constraints \eqref{eq:opf2_rank}. In this regard, we conjecture that the RLS2 is almost always exact in practice. Since we cannot theoretically exclude the possibility of inexactness of OPF in this case, we give some discussions and insights in the next section.

\section{Discussion on Line Constraints }\label{sec:discussion}

\subsection{OPF with Line Constraints}
	Let $\overline{I}_{ij}$ denote the threshold of current through line $i\to j$, then the line constraint can be formulated as
	\begin{equation}\label{eq:current}
	y_{ij}^2(v_i - W_{ij} - W_{ji} + v_j) \le \overline{I}_{ij}^2
	\end{equation}
	which is corresponding to
	\begin{equation}\label{eq:opf1_current}
	y_{ij}^2(V_i - V_j)^2 \le \overline{I}_{ij}^2
	\end{equation}
	in the original problem.
	By adding \eqref{eq:current} into RL1, we have the following model.
	\begin{align*}
	\text{RL2:} ~ \min ~ & h(p) \nonumber \\
	\text{over:} ~ & p,v,W \nonumber \\ 
	\text{s.t.} ~ & \eqref{eq:opf2_p_i}\text{ -- }\eqref{eq:opf2_psd}, \eqref{eq:current}
	\end{align*}
	
	Using the transformation in Theorem \ref{theorem:eq_conv}, the line constraint \eqref{eq:current} can be  be  added into RLS1, yielding the following model.
	\begin{align}
	\text{RLS2:} ~ \min ~ & h(p) \nonumber \\
	\text{over:} ~ & p,P,v,l \nonumber \\
	\text{s.t.} ~ & \eqref{eq:rls1_p_i}\text{ -- }\eqref{eq:rls1_rank} \nonumber \\
	& l_{ij} \le \overline{I}_{ij}^2
	\end{align}

\subsection{Exactness Conditions When Considering Line Constraints}



We first introduce the following theorem to give the conditions that guarantee the exactness of RL2.

\begin{theorem}\label{theorem:rl2_le}
	Assume Assumptions \ref{assumption1} and \ref{assumption2} hold for RL2. Then RL2 is exact under either of the following conditions.
\begin{condition}
   $\forall (i\to j)\in\mathcal{E}$, constraint \eqref{eq:current} is not binding  ;
\end{condition}
\begin{condition}
  For each line $(s\to t)\in\mathcal{E}$, if constraint \eqref{eq:current} is binding for this line, then the corresponding lower bounds on $p_s$ and $p_t$ are not binding.
\end{condition}
\end{theorem}

The proof of Theorem \ref{theorem:rl2_le} can be found in Appendix \ref{proof_thm_rl2_le}.

Furthermore, if the solution to RL2 is not exact, it is possible to find an approximate solution by relaxing the bounds on the power injections, as the following theorem indicates. 

\begin{theorem}\label{theorem:rl2_le_eq}
	Assume Assumptions \ref{assumption1} and \ref{assumption2} hold for RL2 and let $(p,v,W)$ be a feasible solution to RL2. If
	\begin{enumerate}
		\item $(p,v,W)$ violates the rank constraint \eqref{eq:opf2_rank} for a certain line $(s\to t)\in\mathcal{E}$ where the line constraint \eqref{eq:current} is binding;
		\item at least one of $p_s=\underline{p}_s$ and $p_t=\underline{p}_t$ is satisfied;
	\end{enumerate}
	then there must exist  another solution $(p',v,W')$ such that
	\begin{enumerate}
		\item satisfies all the constraints of RL2 except for \eqref{eq:opf2_p_gen} (i.e., \eqref{eq:opf2_p_i}, \eqref{eq:opf2_v} -- \eqref{eq:opf2_psd} and \eqref{eq:current});
		\item satisfies the rank constraint \eqref{eq:opf2_rank}.
		\item has a lower objective value than $(p,v,W)$.
	\end{enumerate}
\end{theorem}

The proof of Theorem \ref{theorem:rl2_le_eq} can be found in Appendix \ref{proof_thm_rl2_le_eq}.


Theorem \ref{theorem:rl2_le_eq} indicates that if the adjustment of power injection $y_{st}\epsilon$ is allowed at bus $s$ and bus $t$, then a solution which satisfies the rank constraint \eqref{eq:opf2_rank} can be constructed. In a DC-MG, such adjustment may be achieved by employing demand response or energy storage.

\subsection{Constructing Approximate Optimal Solutions}
Since RL2 is not always exact, we can check the solution after solving RL2. If the solution satisfies the rank constraint \eqref{eq:opf2_rank}, then it is global optimal for the original problem OPF1 \eqref{eq:opf1}. Otherwise, inspired by the recovery methods of semidefinite program (SDP) relaxation \cite{Fazelnia2017}, we propose two heuristic methods to construct a nearly optimal solution.

Theorem \ref{theorem:rl2_le} indicates that if $(p^*,v^*,W^*)$ is the optimal solution to RL2, and violates the rank constraint \eqref{eq:opf2_rank} for a certain line $(s\to t)\in\mathcal{E}$ where the line constraint \eqref{eq:current} is binding, then at least one of $p_s$ and $p_t$ must reach its lower bound. Otherwise, we can always find another feasible point, which has a smaller objective value than $(p^*,v^*,W^*)$. Therefore, we only need to discuss the cases that at least one of $p_s$ and $p_t$ reaches its lower bound. 

\subsubsection{Direct construction method}
In \cite{Fazelnia2017}, a direct recovery method is proposed for AC networks, using the first column of the optimal solution matrix $W^*$ to recover a nearly optimal solution. It inspires a direct construction method for DC-MGs. 

Assume $(p^*,v^*,W^*)$ is the optimal solution to RL2, which violates the rank constraint \eqref{eq:opf2_rank} for some lines. Instead of using the first column of $W^*$ (see in \cite{Fazelnia2017}), we use $v^*$ to construct an approximate solution $(p,V)$ to the original problem OPF1. 
First, $V_i(i\in\mathcal{N})$ is derived by letting $V_i=\sqrt{v^*_i}$. Then, $p_i(i\in\mathcal{N})$ is derived by substituting $V_i$ into the power balance equation \eqref{eq:opf1_p_i}. Next we show the approximate solution $(p,V)$: 
	\begin{enumerate}
		\item satisfies power balance equation \eqref{eq:opf1_p_i}, the voltage constraint \eqref{eq:opf1_v}, line constraint \eqref{eq:opf1_current};
		\item may violate \eqref{eq:opf1_p_gen} for some $i\in\mathcal{N}$, but the violations have  limited bounds.
	\end{enumerate}

According to the construction of $(p,V)$, it is straightforward to check that $(p,V)$ satisfies \eqref{eq:opf1_p_i} and \eqref{eq:opf1_v}. Hence we only need to examine constraint \eqref{eq:opf1_current} to justify the first assertion.  

For any $i\in\mathcal{N}$, let $\mathcal{C}_i$ denote the set of buses immediately connected to bus $i$. Let $\mathcal{B}_i \subseteq \mathcal{C}_i$  denote the subset such that for any $j\in\mathcal{B}_i$, line $i\sim j$ violates the rank constraint \eqref{eq:opf2_rank}. If all the neighboring lines of bus $i$ do not violate the rank constraint \eqref{eq:opf2_rank}, we have $\mathcal{B}_i=\emptyset$. 
It follows that for all $(i\sim j)\in\mathcal{E}$, $\sqrt{v^*_iv^*_j}>W^*_{ij}$ for $j\in\mathcal{B}_i$, while $\sqrt{v^*_iv^*_j}=W^*_{ij}$ for $j\notin\mathcal{B}_i$. Since $(p^*,v^*,W^*)$ satisfies \eqref{eq:current}, for any $i\sim j(j\in\mathcal{B}_i)$, it follows that
\begin{align*}
	y_{ij}^2(V_i - V_j)^2 &= y_{ij}^2(V_i^2 - 2V_iV_j + V_j^2) \\
	&= y_{ij}^2(v^*_i - 2\sqrt{v^*_iv^*_j} + v^*_j) \\
	&< y_{ij}^2(v^*_i - 2W^*_{ij} + v^*_j) \le \overline{I}_{ij}^2
\end{align*}
and for any $i\sim j(j\notin\mathcal{B}_i)$
\begin{align*}
	y_{ij}^2(V_i - V_j)^2 &= y_{ij}^2(V_i^2 - 2V_iV_j + V_j^2) \\
	&= y_{ij}^2(v^*_i - 2\sqrt{v^*_iv^*_j} + v^*_j) \\
	&= y_{ij}^2(v^*_i - 2W^*_{ij} + v^*_j) \le \overline{I}_{ij}^2.
\end{align*}
Therefore, $(p,V)$ satisfies \eqref{eq:opf1_current}.

Next, we check constraint \eqref{eq:opf1_p_gen}. 
It follows from \eqref{eq:opf1_p_i} that 
\begin{align*}
  p_i &= \sum_{j:j\sim i}V_i\left(V_i - V_j\right)y_{ij} = \sum_{j:j\sim i}\sqrt{v^*_i}\left(\sqrt{v^*_i} - \sqrt{v^*_j}\right)y_{ij} \\
    &= \sum_{j:j\sim i}\left(v^*_i - \sqrt{v^*_iv^*_j}\right)y_{ij} \\
    &= \sum_{j:j\sim i,j\notin\mathcal{B}_i}\left(v^*_i - \sqrt{v^*_iv^*_j}\right)y_{ij} \\
    &+ \sum_{j:j\sim i, j\in\mathcal{B}_i} \left(v^*_i - \sqrt{v^*_iv^*_j}\right)y_{ij} \\
    &\le \sum_{j:j\sim i,j\notin\mathcal{B}_i}\left(v^*_i - W^*_{ij}\right)y_{ij} + \sum_{j:j\sim i, j\in\mathcal{B}_i} \left(v^*_i - W^*_{ij}\right)y_{ij} \\
    &= p^*_i.
\end{align*}
Hence, $p_i$ may violate  \eqref{eq:opf1_p_gen}. And the violation is bounded by
\begin{equation*}
  p_i - p^*_i \le  \sum_{j:j\sim i, j\in\mathcal{B}_i} \left(- \sqrt{v^*_iv^*_j} + W^*_{ij}\right)y_{ij}.
\end{equation*}


\subsubsection{Slack variable method}

\begin{figure}[htbp]
\centering
\includegraphics[width=2.7in]{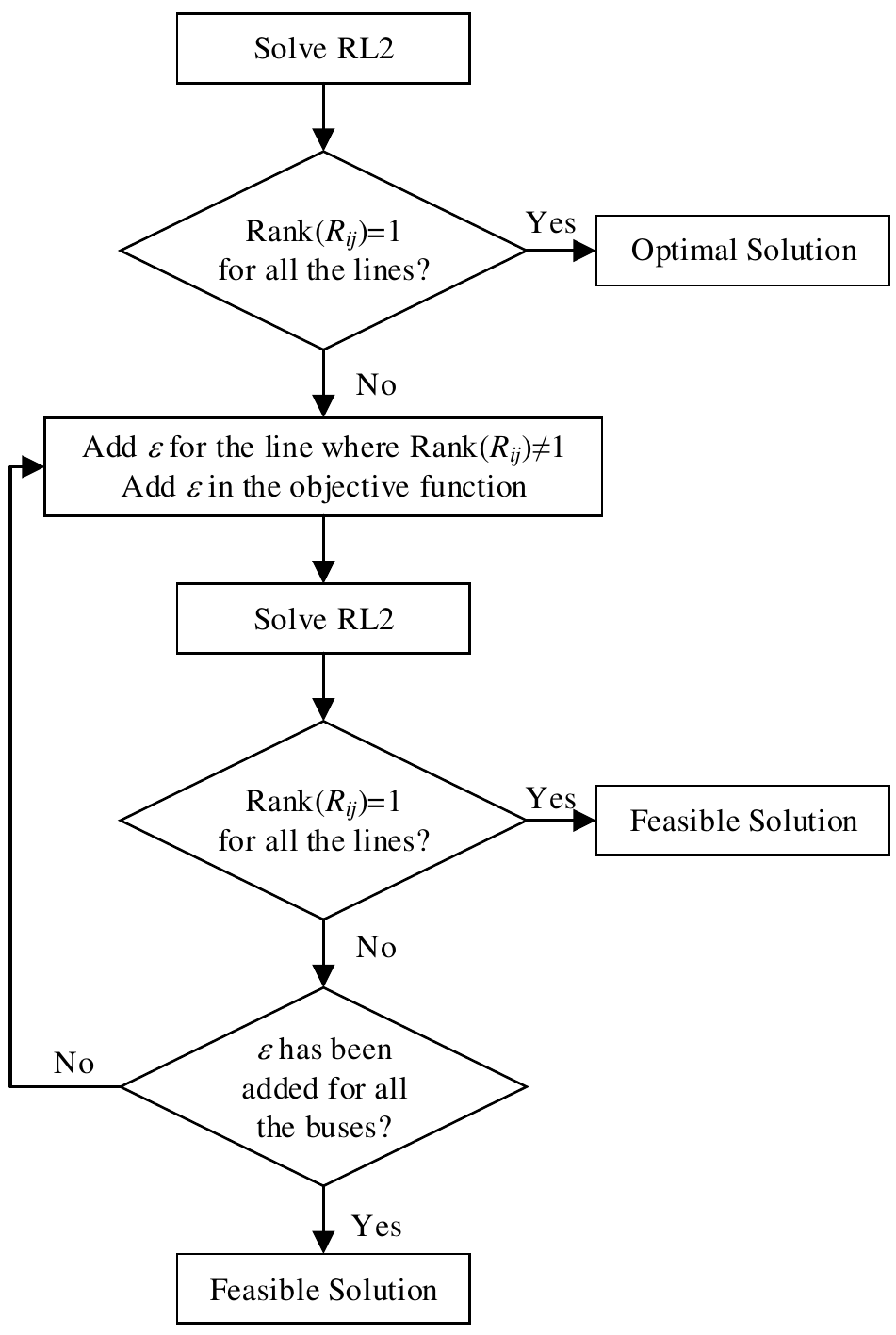}
\caption{Flowchart of slack variable method.}
\label{fig:alg2}
\end{figure}

According to Theorem \ref{theorem:rl2_le}, RL2 is exact if the lower bounds of power injections are not binding. Thus, after solving RL2, if the solution violates the rank constraint \eqref{eq:opf2_rank} for a certain line $(s\sim t)\in\mathcal{E}$, and any lower bound of $p_s$ and $p_t$ is binding, for example, $p_s=\underline{p}_s$, then a corresponding slack variable $\varepsilon_s$ ($\varepsilon_s>0$) can be added into \eqref{eq:opf2_p_gen} to reformulate the constraint as
\begin{equation}
	\underline{p}_s \le p_s - \varepsilon_s \le \overline{p}_s
\end{equation}
so that $p_s$ will not reach its lower bound in the next iteration to solve RL2. Additionally, in order to minimize $\varepsilon_s$, it is also added into the objective function as
\begin{equation*}
\min h(p) + \sum_{i\in\hat{\mathcal{N}}} \varepsilon_i
\end{equation*}
where $\hat{\mathcal{N}}$ is the set of buses where the rank constraint \eqref{eq:opf2_rank} is violated and the power injection lower bound is binding at the same time. 
The procedure is shown in Fig. \ref{fig:alg2}.

\section{Conclusions}\label{sec:conclusion}
In this paper, we have proposed an SOCP relaxation of the OPF problem in stand-alone DC-MGs, which do not consist of any substation with an unconstrained power injection and a fixed nodal voltage. We have also proposed two mild assumptions which only require uniform voltage upper bounds and positive network losses. Under such assumptions, we have proven the exactness of the SOCP relaxation by extending the results in \cite{Gan2014}. Combining the results in \cite{Gan2014} and those in this paper, we have a more comprehensive understanding on the OPF problem in DC-MGs:

\begin{enumerate}
\item 
Under the proposed assumptions, the exactness of SOCP relaxation in DC-MGs does not rely on the operating mode of grid. This property facilitates the optimal operation of DC-MG, since it allows the operator to  achieve optimal operation whether the DC-MG is working in grid-connected or stand-alone mode.
\item 
Unlike AC systems, the exactness of SOCP relaxation in DC-MGs does not rely on the topology of networks. It implies that when the topology is changed due to line switching, the OPF problem can still be converted equivalently into a convex SOCP problem.
\item 
The uniqueness of the global optimal solution to the SOCP problem is independent of the topology and operating mode of DC-MGs. This property is especially useful for setting a target in optimal control of DC-MGs.
\end{enumerate}

\ifCLASSOPTIONcaptionsoff
\newpage
\fi

\bibliographystyle{IEEEtran}
\bibliography{library}

\clearpage
\appendix

\subsection{Proof of Lemma \ref{lemma:trans}}\label{sec:proof_lemma_trans}

\begin{IEEEproof}
\emph{Existence:} Let $V_i=\sqrt{v_i}$ for $i\in \mathcal{N}$. It suffices to show that $V$ satisfies $V_i>0$ for $i\in\mathcal{N}$ and \eqref{eq:slack_var}.

Since $v_i>0$ for $i\in \mathcal{N}$, it follows that $V_i=\sqrt{v_i}>0$ for $i\in \mathcal{N}$. It is straightforward to check that $V$ satisfies \eqref{eq:slack_v}. Since $R_{ij}$ is not full rank, we have
\begin{equation*}
  v_iv_j - W_{ij}W_{ji} = 0, ~ i\to j
\end{equation*}
Since $W_{ij}\ge 0$, we have
\begin{equation*}
W_{ij} = \sqrt{W_{ij}^2} = \sqrt{W_{ij}W_{ji}} = \sqrt{v_iv_j} = V_iV_j
\end{equation*}
for $i\sim j$, i.e., $V,$ satisfies \eqref{eq:slack_w}. This completes the proof of existence.

\emph{Uniqueness:} Let $\tilde{V}$ denote an arbitrary solution to $\tilde{V}_i > 0$ for $i\in \mathcal{N}$ and \eqref{eq:slack_var}. It suffices to show that $\tilde{V}_i=\sqrt{v_i}$ for $i\in \mathcal{N}$.
Assume $\tilde{V}_i\ne\sqrt{v_i}$ for some $i\in \mathcal{N}$, then it follows from \eqref{eq:slack_v} that $\tilde{V}_i=-\sqrt{v_i}<0$, which contradicts with $\tilde{V}_i > 0$ for $i\in \mathcal{N}$. Thus, $\tilde{V}_i=\sqrt{v_i}$ for $i\in \mathcal{N}$. This completes the proof of uniqueness.
\end{IEEEproof}

\subsection{Proof of Theorem \ref{theorem:eq_trans}}
 \label{sec:proof_thm_trans}
\begin{IEEEproof}
	Let $\mathcal{F}_{OPF1}$ and $\mathcal{F}_{OPF2}$ denote the feasible sets of OPF1 and OPF2, respectively.
	Since OPF1 and OPF2 have the same objective function, it suffices to show that there exists a one-to-one map between $\mathcal{F}_{OPF1}$ and $\mathcal{F}_{OPF2}$.
	Specifically,  we show  the map $f:(V)\mapsto (v,W)$ given by \eqref{eq:slack_var} is one-to-one , since $p$ is uniquely determined by $V$ or $(v,W)$. To this end, it suffices to show  the map is both into and onto. On the one hand, it is straightforward that for any $V_1\ne V_2\in\mathcal{F}_{OPF1}$, $f(V_1)\ne f(V_2)\in\mathcal{F}_{OPF2}$. On the other hand, as Lemma \ref{lemma:trans} says, for each $(v,W)\in\mathcal{F}_{OPF2}$, there exists a unique $f^{-1}(v,W)=V\in\mathcal{F}_{OPF1}$. This completes the proof.
\end{IEEEproof}

\subsection{Proof of Lemma \ref{lemma:v}}\label{sec:proof_lemma_v}

\begin{IEEEproof}
For brevity, we only prove the case of bus $s$ as the proof of bus $t$ is the same. 

With regard to the value of $\underline{p}_s$, we only need to discuss two cases: $p_s=\underline{p}_s<0$ and $p_s=\underline{p}_s=0$.

In the first case, we have
\begin{equation*}
  p_s = \sum_{i:i\sim s}\left(v_s - W_{si}\right)y_{si} < 0
\end{equation*}
according to \eqref{eq:opf2_p_i}. Therefore, $\left(v_s - W_{si}\right)y_{si}< 0$ for some $i^{\flat}\in \mathcal{N}$. It follows from \eqref{eq:opf2_w_0} -- \eqref{eq:opf2_psd} that $W_{si^{\flat}}\le\sqrt{v_sv_{i^{\flat}}}$. According to Assumption \ref{assumption1}, $\overline{V}_s=\overline{V}_{i^{\flat}}$. Thus, we have
\begin{equation*}
  v_s < W_{si^{\flat}} \le \sqrt{v_sv_{i^{\flat}}} \le \sqrt{\overline{V_s}^2\overline{V}_{i^{\flat}}^2} = \overline{V}_s^2.
\end{equation*}

In the second case, we have 
\begin{equation}
\label{eq:powerblance}
  p_s = \sum_{i:i\sim s}\left(v_s - W_{si}\right)y_{si} = 0
\end{equation}
according to \eqref{eq:opf2_p_i}. Then we discuss the following two cases: 

1) When $v_s-W_{si}=0$ for all $i\sim s$, we have $v_s=W_{st}$ for $i=t$. Since $(p,v,W)$ violates the rank constraint \eqref{eq:opf2_rank} for $s\to t$, $R_{st}$ is non-singular. Hence, $W_{st}\ne\sqrt{v_sv_t}\Rightarrow W_{st}<\sqrt{v_sv_t}$ due to \eqref{eq:opf2_w_0} -- \eqref{eq:opf2_psd}. According to Assumption \ref{assumption1}, $\overline{V}_s=\overline{V}_t$. It follows that
$v_s = W_{st} < \sqrt{v_sv_t} \le \sqrt{\overline{V_s}^2\overline{V}_t^2} = \overline{V}_s^2.$

2) When $v_s-W_{si}=0$ does not hold for some of $i\sim s$, there must exist at least one $i^{\flat}\sim s$ such that $v_s-W_{si^{\flat}}<0$ (and accordingly there exists at least another $i^{\sharp}\sim s$ such that $v_s-W_{si^{\sharp}}>0$ in order to satisfy \eqref{eq:powerblance} ). Hence we have
$v_s < W_{si^{\flat}} \le \sqrt{v_sv_{i^{\flat}}} \le \sqrt{\overline{V_s}^2\overline{V}_{i^{\flat}}^2} = \overline{V}_s^2$ due to $\overline{V}_s=\overline{V}_{i^{\flat}}$.
This completes the proof.
\end{IEEEproof}

\subsection{Proof of Lemma \ref{lemma:le}}\label{sec:proof_lemma_le}

\begin{IEEEproof}
Since $(p,v,W)$ satisfies \eqref{eq:opf2_w_0} -- \eqref{eq:opf2_psd}, we have $0\le W_{ij}\le\sqrt{v_iv_j}$ for $i\to j$. Since $(p,v,W)$ violates the rank constraint \eqref{eq:opf2_rank} for $s\to t$, we have $W_{st}\ne\sqrt{v_sv_t}$. Thus, $W_{st}<\sqrt{v_sv_t}$. We can always choose a small enough number $\epsilon>0$ such that
\begin{equation*}
\epsilon < \min\left\{\frac{p_s-\underline{p}_s}{y_{st}}, \frac{p_t-\underline{p}_t}{y_{st}}, \sqrt{v_sv_t}-W_{st}\right\}.
\end{equation*}

Following the line of the proof of Lemma 7 in \cite{Gan2014}, we can use $\epsilon$ to construct $W'$ as
\begin{equation*}
  W'_{ij} = \begin{cases}
  W_{ij} + \epsilon & \text{if} ~ \{i,j\}=\{s,t\}; \\
  W_{ij} & \text{otherwise};
  \end{cases}
\end{equation*}
and construct $p'$ as
\begin{equation*}
  p'_i = \sum_{j:j\sim i}\left(v_i - W'_{ij}\right)y_{ij},~ i\in \mathcal{N}.
\end{equation*}

The point $(p',v,W')$ satisfies \eqref{eq:opf2_p_i} due to the construction of $p'$. When $i\ne s,t$, we have
\begin{equation*}
  p'_i = \sum_{j:j\sim i}\left(v_i - W'_{ij}\right)y_{ij} = \sum_{j:j\sim i}\left(v_i - W_{ij}\right)y_{ij} = p_i.
\end{equation*}
When $i=s,t$, we have 
\begin{align*}
  p'_i &= \sum_{j:j\sim i}(v_i - W'_{ij})y_{ij} \\
  &= \sum_{j:j\sim i}(v_i - W_{ij})y_{ij} - y_{st}\epsilon \\
  &= p_i - y_{st}\epsilon \in (\underline{p}_i,\overline{p}_i).
\end{align*}
Hence, $(p',v,W')$ satisfies \eqref{eq:opf2_p_gen} and
\begin{equation*}
  p'_i \begin{cases}
  < p_i & \text{if}\quad i = s,t; \\
  = p_i & \text{otherwise}.
  \end{cases}
\end{equation*}
It follows that $h(p')<h(p)$ since $f$ is strictly increasing in $p_i$ for $i\in\mathcal{N}$. The point $(p',v,W')$ satisfies \eqref{eq:opf2_v} since $v$ remains the same as the feasible point $(p,v,W)$. The point $(p',v,W')$ satisfies \eqref{eq:opf2_w_0} due to the construction of $W'$. The point $(p',v,W')$ satisfies \eqref{eq:opf2_w} since
\begin{equation*}
  W'_{ij} - W'_{ji} = W_{ij} - W_{ji} = 0
\end{equation*}
for $i\to j$. 
Additionally, the point $(p',v,W')$ satisfies \eqref{eq:opf2_psd} since
\begin{equation*}
  W'_{ij} = W_{ij} \in [0,\sqrt{v_iv_j}]
\end{equation*}
when $\{i,j\}\ne\{s,t\}$ and
\begin{equation*}
  W'_{ij} = W_{ij} + \epsilon \in [0,\sqrt{v_iv_j})
\end{equation*}
when $\{i,j\}=\{s,t\}$. This completes the proof.
\end{IEEEproof}

\subsection{Proof of Lemma \ref{lemma:le_eq}}\label{sec:proof_lemma_le_eq}

\begin{IEEEproof}
We present the proof for the case ($p_s=\underline{p}_s$ \& $p_t>\underline{p}_t$). The proof for the case ($p_s>\underline{p}_s$ \& $p_t=\underline{p}_t$) is similar and omitted for brevity.

Similar to the proof of Lemma \ref{lemma:le}, we have $W_{st}<\sqrt{v_sv_t}$. Additionally, it follows from Lemma \ref{lemma:v} that $v_s<\overline{V}_s^2$, so we can always choose a small enough number $\epsilon>0$ such that 
\begin{equation*}
  \epsilon < \min \left\{\frac{p_t-\underline{p}_t}{y_{st}}, \sqrt{v_sv_t}-W_{st}, \frac{\sum_{j:j\sim s}y_{sj}}{y_{st}}\left(\overline{V}_s^2 - v_s\right) \right\},
\end{equation*}
then
\begin{align*}
& p_t - y_{st}\epsilon > \underline{p}_t, \\
& W_{st} + \epsilon < \sqrt{v_sv_t}, \\
& v_s + \frac{y_{st}}{\sum_{j:j\sim s}y_{sj}}\epsilon < \overline{V}_s^2.
\end{align*}

Following the line of the proof of Lemma 8 in \cite{Gan2014}, we can use $\epsilon$ to construct $W'$ as
\begin{equation*}
  W'_{ij} = \begin{cases}
  W_{ij} + \epsilon & \text{if}~\{i,j\} = \{s,t\}; \\
  W_{ij} & \text{otherwise};
  \end{cases}
\end{equation*}
construct $v'$ as
\begin{equation*}
  v'_{i} = \begin{cases}
  v_{i} + \frac{y_{st}}{\sum_{j:j\sim i}y_{ij}}\epsilon & \text{if}\quad i = s; \\
  v_{i} & \text{otherwise};
  \end{cases}
\end{equation*}
and construct $p'$ as
\begin{equation*}
  p'_i = \sum_{j:j\sim i}(v'_i - W'_{ij})y_{ij},~ i\in \mathcal{N}.
\end{equation*}

The point $(p',v',W')$ satisfies \eqref{eq:opf2_p_i} according to the construction of $p'$. When $i\ne s,t$, we have
\begin{equation*}
  p'_i = \sum_{j:j\sim i}\left(v'_i - W'_{ij}\right)y_{ij}
  = \sum_{j:j\sim i}\left(v_i - W_{ij}\right)y_{ij} = p_i.
\end{equation*}
When $i=s$, we have
\begin{align*}
  p'_s &= \sum_{j:j\sim s}\left(v'_s - W'_{sj}\right)y_{sj} \\
  &= \sum_{j:j\sim s}\left(v_s - W_{sj}\right)y_{sj} + (v'_s-v_s)\sum_{j:j\sim s}y_{sj} - y_{st}\epsilon\\
  &= \sum_{j:j\sim s}\left(v_s - W_{sj}\right)y_{sj} = p_s.
\end{align*}
When $i=t$, we have
\begin{align*}
  p'_t &= \sum_{j:j\sim t}(v'_t - W'_{tj})y_{tj} \\
  &= \sum_{j:j\sim t}(v_t - W_{tj})y_{tj} - y_{st}\epsilon \\
  &= p_t - y_{st}\epsilon \in (\underline{p}_t,\overline{p}_t).
\end{align*}
Thus, $(p',v',W')$ satisfies \eqref{eq:opf2_p_gen} and 
\begin{equation*}
  p'_i \begin{cases}
  < p_i & \text{if}\quad i = t; \\
  = p_i & \text{otherwise}.
  \end{cases}
\end{equation*}
It follows that $h(p')<h(p)$. Since 
$v'_i=v_i$ if $i\ne s$ and $v_i<v_i'<\overline{V}_i^2$ if $i=s$, the point $(p',v',W')$ satisfies \eqref{eq:opf2_v}. 

Similar to the proof of Lemma \ref{lemma:le}, it is straightforward to check that the point $(p',v',W')$ satisfies \eqref{eq:opf2_w_0}, \eqref{eq:opf2_w} and \eqref{eq:opf2_psd}.
This completes the proof.
\end{IEEEproof}

\subsection{Proof of Lemma \ref{lemma:eq}}\label{sec:proof_lemma_eq}

\begin{IEEEproof}
Similar to the proof of Lemma \ref{lemma:le}, we have $W_{st}<\sqrt{v_sv_t}$. According to Lemma \ref{lemma:v}, we have
\begin{equation*}
  v_s<\overline{V}_s^2,v_t<\overline{V}_t^2.
\end{equation*}
Therefore, we can always choose a small enough number $\epsilon>0$ so that
\begin{align*}
  \epsilon < \min & \left\{\sqrt{v_sv_t}-W_{st}, \frac{\sum_{j:j\sim s}y_{sj}}{y_{st}}\left(\overline{V}_s^2 - v_s\right), \right. \\
  & \left.\frac{\sum_{j:j\sim t}y_{tj}}{y_{st}}\left(\overline{V}_t^2 - v_t\right) \right\},
\end{align*}
then
\begin{align*}
& W_{st} + \epsilon < \sqrt{v_sv_t}, \\
& v_s + \frac{y_{st}}{\sum_{j:j\sim s}y_{sj}}\epsilon < \overline{V}_s^2, \\
& v_t + \frac{y_{st}}{\sum_{j:j\sim t}y_{tj}}\epsilon < \overline{V}_t^2.
\end{align*}

Following the line of the proof of Lemma 9 in \cite{Gan2014}, we can use $\epsilon$ to construct $W'$ as
\begin{equation*}
  W'_{ij} = \begin{cases}
  W_{ij} + \epsilon & \text{if} ~ \{i,j\}=\{s,t\}; \\
  W_{ij} & \text{otherwise};
  \end{cases}
\end{equation*}
and construct $v'$ as
\begin{equation} \label{eq:lemma2_v}
  v'_i = \begin{cases}
  v_i + \frac{y_{st}}{\sum_{j:j\sim i}y_{ij}}\epsilon & i = s,t; \\
  v_i & \text{otherwise}.
  \end{cases}
\end{equation}

The point $(p,v',W')$ satisfies \eqref{eq:opf2_p_i} since when $i=s,t$,
\begin{align*}
  & \sum_{j:j\sim i}\left(v'_i - W'_{ij}\right)y_{ij} \\
  &= \sum_{j:j\sim i}\left(v_i - W_{ij}\right)y_{ij} + \sum_{j:j\sim i}\frac{y_{st}}{\sum_{j:j\sim i}y_{ij}}\epsilon y_{ij} - \epsilon y_{st} \\
  &= \sum_{j:j\sim i}\left(v_i - W_{ij}\right)y_{ij} = p_i.
\end{align*}
When $i\ne s,t$, we have
\begin{equation*}
  \sum_{j:j\sim i}\left(v'_i - W'_{ij}\right)y_{ij} = \sum_{j:j\sim i}\left(v_i - W_{ij}\right)y_{ij} = p_i.
\end{equation*}
Since $p$ does not change, the point $(p,v',W')$ satisfies \eqref{eq:opf2_p_gen}. 
It follows from \eqref{eq:lemma2_v} that $v'_i=v_i$ if $i\ne s,t$ and $v_i<v_i'<\overline{V}_i^2$ if $i=s,t$. Hence, the point $(p,v',W')$ satisfies \eqref{eq:opf2_v}. 
Similar to the proof of Lemma \ref{lemma:le}, it is straightforward to check that the point $(p',v',W')$ satisfies \eqref{eq:opf2_w_0}, \eqref{eq:opf2_w} and \eqref{eq:opf2_psd}.
Since $\mathcal{G}$ is connected, there always exists some $i \to j$ such that $\{i,j\}\ne\{s,t\}$ and $\{\{i\},\{j\}\}\cap\{\{s\},\{t\}\}\ne\emptyset$. Then we have 
\begin{equation*}
  W'_{ij}=W_{ij}\le\sqrt{v_iv_j}<\sqrt{v_i'v_j'}.
\end{equation*}
Particularly, when $\{i,j\}=\{s,t\}$, $W'_{ij}=W_{ij}+\epsilon<\sqrt{v_iv_j}<\sqrt{v_i'v_j'}$. It means $(p,v',W')$ violates the rank constraint \eqref{eq:opf2_rank} for $s\to t$ and all its neighboring lines. 
This completes the proof.
\end{IEEEproof}

\subsection{Proof of Theorem \ref{theorem:exact_relax}}
\label{proof_thm_relax}
\begin{IEEEproof}
	To prove that RL1 is exact, it suffices to show that any optimal solution to RL2 satisfies the rank constraint \eqref{eq:opf2_rank}. We first assume RL1 is \textit{not} exact for the sake of contradiction. Then there must exist at least one optimal solution $(p^{*},v^{*},W^{*})$ of RL, which violates the rank constraint \eqref{eq:opf2_rank} for a certain line $(s\to t)\in\mathcal{E}$. We discuss from the following three aspects of the power injections $p_s$ and $p_t$ corresponding to $s\to t$. 
	
	(1)  $p^{*}_{s}>\underline{p}_{s}$ \& $p^{*}_{t}>\underline{p}_{t}$. In this case, as Lemma \ref{lemma:le} points out, there must exist a feasible $(p',v',W')$ which has a smaller objective value than $(p^{*},v^{*},W^{*})$. In this situation, $(p^{*},v^{*},W^{*})$ cannot be optimal for RL. 
	
	(2) ($p^{*}_s=\underline{p}_s$ \& $p^{*}_t>\underline{p}_t$) or ($p^{*}_s>\underline{p}_s$ \& $p^{*}_t=\underline{p}_t$). In this case, Lemma \ref{lemma:le_eq} states that there must exist a feasible $(p',v',W')$ that has a smaller objective value than $(p^{*},v^{*},W^{*})$. Hence, $(p^{*},v^{*},W^{*})$ cannot be optimal for RL. 
	
	(3) $p^{*}_s=\underline{p}_s$ \&  $p^{*}_t=\underline{p}_t$. In this case, we show such an optimal solution $(p^{*},v^{*},W^{*})$ does not exist. According to Lemma \ref{lemma:eq}, there always exists a feasible $(p',v',W')$ which violates the rank constraint \eqref{eq:opf2_rank} for $s\to t$ and all its neighboring lines. Since $\mathcal{G}$ is connected and $p'=p^*$, we can further propagate such construction to obtain another feasible $(p^{*},v^{\dag},W^{\dag})$ of RL, which violates the rank constraint \eqref{eq:opf2_rank} for all the lines. Hence, $p^{*}$ satisfies $p^{*}_i=\underline{p}_i$ for all $i\in \mathcal{N}$. It follows that 
	$$\sum_{i\in \mathcal{N}}p^{*}_i = \sum_{i\in \mathcal{N}}\underline{p}_i \le 0$$
	This contradicts Assumption \ref{assumption2} which states $\sum_{i\in \mathcal{N}}p^{*}_i>0$.  
	Thus, $(p^{*},v^{*},W^{*})$ does not exist. 
	
	Thus, every optimal solution must satisfy the rank constraint \eqref{eq:opf2_rank}, i.e., RL1 is exact. This completes the proof.
\end{IEEEproof}

\subsection{Proof of Theorem \ref{theorem:RL1_unique}}
\label{proof_thm_unique}
\begin{IEEEproof}
	Let $x^{1*}:=(p^{1*},v^{1*},W^{1*})$ and $x^{2*}:=(p^{2*},v^{2*},W^{2*})$ be two optimal solutions of RL1, then we have
	\begin{equation}\label{eq:f_eq}
	\sum_{i\in\mathcal{N}}f_i(p^{1*}_i) = \sum_{i\in\mathcal{N}}f_i(p^{2*}_i)
	\end{equation}
	Following the line of the proof of Theorem 3 in \cite{Gan2014}, we know
	\begin{align*}
	& \frac{v^{1*}_i}{v^{2*}_i} = \frac{v^{1*}_j}{v^{2*}_j} = \eta, ~ i\sim j; \\
	& W^{1*}_{ij} = \sqrt{v^{1*}_iv^{1*}_j} = \eta\sqrt{v^{2*}_iv^{2*}_j} = \eta W^{2*}_{ij}.
	\end{align*}
	It follows from \eqref{eq:opf2_p_i} that
	\begin{align*}
	p^{1*}_i &= \sum_{j:j\sim i}\left(v^{1*}_i - W^{1*}_{ij}\right)y_{ij} \\
	&= \sum_{j:j\sim i}\left(\eta v^{2*}_i - \eta W^{2*}_{ij}\right)y_{ij} \\
	&= \eta p^{2*}_i.
	\end{align*}
	Since $f_i(p_i)$ is strictly increasing, we must have $\eta=1$. Otherwise, it contradicts \eqref{eq:f_eq}. Thus, $x^{1*}=x^{2*}$, i.e., RL1 has at most one optimal solution. This completes the proof.
\end{IEEEproof}

\subsection{Proof of Theorem \ref{theorem:eq_conv}}
\label{proof_thm_eq_conv}
\begin{IEEEproof}
	Let $\mathcal{F}_{RL1}$ and $\mathcal{F}_{RLS1}$ denote the feasible sets of RL1 and RLS1, respectively. 
	It suffices to show that the map $g$ \eqref{eq:convert} is one-to-one between $\mathcal{F}_{RL1}$ and $\mathcal{F}_{RLS1}$, since $p$ is determined by $(v,W)$ in $\mathcal{F}_{RL1}$ and determined by $P$ in $\mathcal{F}_{RLS1}$. 
	On the one hand, it is straightforward that for any $(v_1,W_1)\ne (v_2,W_2)\in\mathcal{F}_{RL1}$, $g(v_1,W_1)\ne g(v_2,W_2)\in\mathcal{F}_{RLS1}$. On the other hand, for any $(v',P,l)\in\mathcal{F}_{RLS1}$, it can be verified that there exists a $g^{-1}(v',P,l)=(v,W)\in\mathcal{F}_{RL1}$. This completes the proof.
\end{IEEEproof}

\subsection{Proof of Theorem \ref{theorem:rl2_le}}
\label{proof_thm_rl2_le}
\begin{IEEEproof}
	When the line constraint \eqref{eq:current} is not binding, RL2 is equivalent to RL1, which is exact under Assumptions \ref{assumption1} and \ref{assumption2} according to Theorem \ref{theorem:exact_relax}. However, when constraint \eqref{eq:current} is binding, Theorem \ref{theorem:exact_relax} does not hold. In this situation, assume $(p^*,v^*,W^*)$ is the optimal solution to RL2. Assume $(p^*,v^*,W^*)$ violates the rank constraint \eqref{eq:opf2_rank} for a certain line $(s\to t)\in\mathcal{E}$ where constraint \eqref{eq:current} is binding, and the lower bounds of $p_s$ and $p_t$ are not binding. 
	Since $(p^*,v^*,W^*)$ satisfies \eqref{eq:opf2_w_0} -- \eqref{eq:opf2_psd}, we have $0\le W^*_{ij}\le\sqrt{v^*_iv^*_j}$ for $i\to j$. Since $(p^*,v^*,W^*)$ violates the rank constraint \eqref{eq:opf2_rank} for $s\to t$, we have $W^*_{st}\ne\sqrt{v^*_sv^*_t}$. Thus, $W^*_{st}<\sqrt{v^*_sv^*_t}$. We can always choose a small enough number $\epsilon>0$ such that

	\begin{equation*}
	\epsilon < \min\left\{\frac{p^*_s-\underline{p}_s}{y_{st}}, \frac{p^*_t-\underline{p}_t}{y_{st}}, \sqrt{v^*_sv^*_t}-W^*_{st}\right\}.
	\end{equation*}
	Then we can use $\epsilon$ to construct $(p',v^*,W')$ where
	\begin{equation*}
	W'_{ij} := \begin{cases}
	W^*_{ij} + \epsilon & \text{if} ~ \{i,j\}=\{s,t\}; \\
	W^*_{ij} & \text{otherwise};
	\end{cases}
	\end{equation*}
	\begin{equation*}
	p'_i := \sum_{j:j\sim i}\left(v^*_i - W'_{ij}\right)y_{ij},~ i\in \mathcal{N}.
	\end{equation*}
	
	It is easy to verify that $(p',v^*,W')$ satisfies \eqref{eq:opf2_p_i} -- \eqref{eq:opf2_psd} and $h(p')<h(p^*)$. Additionally, when $\{i,j\}\ne\{s,t\}$, we have
	\begin{equation*}
	y_{ij}^2(v^*_i - W'_{ij} - W'_{ji} + v^*_j) = y_{ij}^2(v^*_i - W^*_{ij} - W^*_{ji} + v^*_j) \le \overline{I}_{ij}^2.
	\end{equation*}
	When $\{i,j\}=\{s,t\}$, since the line constraint \eqref{eq:current} is binding,
	\begin{equation*}
	y_{ij}^2(v^*_i - W'_{ij} - W'_{ji} + v^*_j) = y_{ij}^2(v^*_i - W^*_{ij} - W^*_{ji} + v^*_j - 2\epsilon) < \overline{I}_{ij}^2.
	\end{equation*}
	Thus, the point $(p',v^*,W')$ satisfies \eqref{eq:current}. It means that the point $(p',v^*,W')$ is feasible for RL2 and has a smaller objective value than $(p^*,v^*,W^*)$. Thus, $(p^*,v^*,W^*)$ cannot be the optimal solution. In contrast, the optimal solution to RL2 must satisfy the rank constraint \eqref{eq:opf2_rank}.
	This completes the proof.
\end{IEEEproof}

\subsection{Proof of Theorem \ref{theorem:rl2_le_eq}}
\label{proof_thm_rl2_le_eq}
\begin{IEEEproof}
	Noting that $(p,v,W)$ satisfies \eqref{eq:opf2_w_0} -- \eqref{eq:opf2_psd}, we have $0\le W_{ij}\le\sqrt{v_iv_j}$ for $i\to j$. Since $(p,v,W)$ violates the rank constraint \eqref{eq:opf2_rank} for $s\to t$, $W_{st}\ne\sqrt{v_sv_t}$. Hence $W_{st}<\sqrt{v_sv_t}$. Then  we can always choose 
	$$	\epsilon := \sqrt{v_sv_t}-W_{st}>0$$ 
	to construct another solution $(p',v,W')$ as below: 
		\begin{equation*}
		W'_{ij} := \begin{cases}
		W_{ij} + \epsilon & \text{if} ~ \{i,j\}=\{s,t\}; \\
		W_{ij} & \text{otherwise};
		\end{cases}
		\end{equation*}
		\begin{equation*}
		p'_i := \sum_{j:j\sim i}\left(v_i - W'_{ij}\right)y_{ij},~ i\in \mathcal{N}.
		\end{equation*}
	
	It is easy to verify that  $(p',v,W')$ satisfies \eqref{eq:opf2_p_i}, \eqref{eq:opf2_v} -- \eqref{eq:opf2_psd} and \eqref{eq:current}. Hence the first assertion is proven.
	
	In terms of the second assertion,  $(p',v,W')$ satisfies \eqref{eq:opf2_rank} by noting  
	\begin{equation*}
	W'_{ij} = W_{ij} = \sqrt{v_iv_j} ~~~~\text{for}~ \{i,j\}\ne\{s,t\}
	\end{equation*}
	and 
	\begin{equation*}
	W'_{ij} = W_{ij} + \epsilon = \sqrt{v_iv_j} ~~~~\text{for}~\{i,j\}=\{s,t\}.
	\end{equation*}
	
	Next we consider the third assertion. 
		According to the construction of $(p',v,W')$, for any $i\ne s,t$, $p'$ satisfies
		\begin{equation*}
		p'_i = \sum_{j:j\sim i}\left(v_i - W'_{ij}\right)y_{ij} = \sum_{j:j\sim i}\left(v_i - W_{ij}\right)y_{ij} = p_i
		\end{equation*}
		and for any when $i=s,t$, 
		\begin{align*}
		p'_i &= \sum_{j:j\sim i}(v_i - W'_{ij})y_{ij} \\
		&= \sum_{j:j\sim i}(v_i - W_{ij})y_{ij} - y_{st}\epsilon \\
		&= p_i - y_{st}\epsilon < p_i
		\end{align*}
		Then $(p',v,W')$ has a lower objective value  than  $(p,v,W)$ because the objective function is strictly increasing in $p_i$.
	This completes the proof.
\end{IEEEproof}

\end{document}